\documentclass[envcountsame]{mmnp}
\usepackage{amsmath,amssymb,graphicx,subfigure}
\usepackage{amsfonts}
\usepackage{graphicx}
\usepackage[latin1]{inputenc}
\usepackage{setspace}
\usepackage{bm}
\usepackage{color}

\numberwithin{equation}{section}

\newcommand{\Real}{\mathbb{R}}

\begin{document}

\title{Leaders do not look back, or do they?}

\author{A. N. Gorban \inst{1}, N. Jarman \inst{1,2}, E. Steur \inst{2} \thanks{E. Steur is now with Eindhoven University of Technology, Institute for Complex Molecular Systems and Department of Mechanical Engineering, the Netherlands}, C. van Leeuwen\inst{2}, I. Yu. Tyukin \inst{1,3} \thanks{\email {I.Tyukin@le.ac.uk}}}

\vspace{0.5cm}

\institute{\inst{1} {University of Leicester, Department of Mathematics, United Kingdom}\\
\inst{2} {KU Leuven, Department of Psychology,\\ Laboratory for Perceptual Dynamics, Belgium}\\
\inst{3} {Saint-Petersburg State Electrotechnical University,\\ Department of Automation and Control Processes, Russia}}

\abstract{We study the effect of adding to a directed chain of interconnected systems a directed feedback from the last element in the chain to the first. The problem is closely related to the fundamental question of how a change in network topology may influence the behavior of coupled systems.
 We begin the analysis by investigating a simple linear system. The matrix that specifies the system dynamics is the transpose of the network Laplacian matrix, which codes the connectivity of the network. Our analysis shows that for any nonzero complex eigenvalue $\lambda$ of this matrix, the following inequality holds: $\frac{|\Im \lambda |}{|\Re \lambda |} \leq
\cot\frac{\pi}{n}$. This bound is sharp, as it becomes an equality for an eigenvalue of a
simple directed cycle with uniform interaction weights. The latter has the slowest decay of oscillations among all other network configurations with the
same number of states. The result is generalized to directed rings and chains of identical nonlinear oscillators. For directed rings, a lower bound $\sigma_c$ for the connection strengths that guarantees asymptotic synchronization is found to follow a similar pattern: $\sigma_c=\frac{1}{1-\cos\left( 2\pi /n\right)} $. Numerical analysis revealed that, depending on the network size $n$, multiple dynamic regimes co-exist in the state space of the system. In addition to the fully synchronous state a rotating wave solution occurs. The effect is observed in networks exceeding a certain critical size. The emergence of a rotating wave highlights the importance of long chains and loops in networks of oscillators: the larger the size of chains and loops, the more sensitive the network dynamics becomes to removal or addition of a single connection.
}

\keywords{coupled systems\sep reaction network\sep eigenvalue\sep synchronization\sep wave solutions}

\subjclass{34A30\sep 34D06\sep 34D45\sep 92B20\sep 92B25}

\titlerunning{Leaders do not look back}

\maketitle

\section{Introduction}

A fundamental question in complex networks is how topology influences the overall network behavior. This issue is crucial for understanding a range of phenomena in elementary chemical kinetic systems, populations of agents, and processes in the neuronal circuits of the human brain \cite{Gaiteri2011}. It is well known that a sufficiently strong diffusive coupling will lead to globally asymptotically stable synchronization in a large class of systems \cite{Pogromsky1998}. Some network topologies, moreover, may give rise to partial synchronization \cite{Pogromsky2002,Belykh2000,Belykh2004}, and networks with dynamically changing topologies were shown to exhibit complex multi-stable dynamics (see e.g. \cite{Kurths2014} and references therein).

Even when network topology is not changing dynamically, its influence on the overall network dynamics is well documented. Examples of how the network topology may affect e.g. coherence of network dynamics are provided in \cite{Belykh2004}. The authors showed that shortcuts in otherwise regular lattices significantly reduce the critical coupling strength needed for achieving global asymptotic stability of the synchronous state. Hence networks with shortcuts can be considered as more efficient than regular ones in term of resources spent, such as the total number of connections and their strength, for reaching and maintaining synchronous regimes. This aspect of shortcuts appear to be crucial for forming small-world structures \cite{Gong2004,Jarman2014} in evolving networks. Further examples showing significant dependence of network dynamics on the corresponding connectivity graphs can be found in  \cite{Gaiteri2011,Marttunen}.


Understanding the problem of how the network topology affects its dynamics is a huge theoretical and practical challenge if considered in its full generality. Here we will focus on a much simpler objective. In particular we will discuss and analyze two basic and extreme topologies which any network will contain as a subgraph: a directed chain, and a directed cycle. Not only may these be considered as basic building blocks of arbitrary network topologies; many networks can be reduced to chains and cycles as well \cite{Gorban2008,Gorban2010}. Moreover, recent computational studies revealed that cycles could be important on their own for sustaining coherent oscillatory network activity \cite{Garcia2014}.

We begin our investigation by analysing the dynamics of a system of coupled neutrally stable linear equations. The dynamics are essentially governed by a coupling matrix that corresponds to directed interconnections in the system. The equations can also be viewed as a model describing the dynamics of damped oscillations in kinetic systems, and thus in what follows we refer to it as such. The results are provided in Section \ref{sec:eigenvalues}. In Section \ref{sec:neural} we consider a more generalized setting, in which the dynamics of each individual node is governed by a nonlinear, albeit semi-passive \cite{Pogromsky1998}, oscillator. The equations describing oscillators in each node are of the FitzHugh-Nagumo (FHN) type \cite{FitzHugh}. These oscillators are, in turn, an adaptation of the van der Pol oscillator \cite{van_der_Pol}. We found that for the case of the directed cycle the value of critical coupling needed to maintain globally asymptotically  stable synchrony is $O\left(\frac{1}{1-\cos\left( 2\pi /n\right)}\right)$, whereas the synchronization threshold for systems organized into directed chains does not depend on $n$. Moreover, the error dynamics corresponding to the simple directed cycle rapidly becomes underdamped for large $n$; this enables resonances between the dynamics of individual nodes and the coupling dynamics. Further numerical analysis show that not only fully synchronous oscillations may occur in these types of network, but also a stable rotating wave solution may emerge. These two dynamic regimes co-exist for a broad range of coupling strength and for number of systems. Occasionally we observe a shift towards prevalence of rotating waves for large enough $n$.   Section \ref{sec:discussion} contains a discussion of our findings, and  Section \ref{sec:conclusion} concludes the paper.


\setcounter{equation}{0}

\section{Coupled Neutrally Stable Systems}\label{sec:eigenvalues}
Consider the following system of linear first-order differential equations:
\begin{equation}\label{eq:sysP}
\dot P = K P,
\end{equation}
 where $P = \mathrm{col}(p_1,p_2,\ldots,p_n) \in \mathbb{R}^n$, and matrix $K=(k_{ij})$ is defined as  follows:
\begin{equation}\label{eq:K}
k_{ij}=\left\{\begin{array}{cl}
q_{ij}, \ q_{ij}\geq 0 &\mbox{ if } i\neq j;\\
-\sum_{m, \, m\neq i} q_{m i} &\mbox{ if } i = j .
\end{array}\right.
\end{equation}
Note that $K$ is a Metzler matrix\footnote{A Metzler matrix is a matrix with non-negative off-diagonal entries} with zero column sums. Off-diagonal elements $k_{ij}$, $i\neq j$, of the matrix $K$ can be viewed as the connection weights between the $i$-th and the $j$-th nodes in the network. The matrix $K$ can be related to the Laplacian matrix $L$ (see e.g. \cite{Bollobas}) of an associated directed network in which the overall connectivity pattern is the same except for that the direction of all connections is altered. The Laplacian for the latter network is thus $L=-K^T$. Note, however, that this relation does not necessarily hold for the original network.

\textcolor{black}{System \eqref{eq:sysP} is a commonly used model of first-order kinetics with a finite number of states. In this case the variables $p_i$ may represent {\em concentration},  {\em probability}, or {\em population} of these states.} \textcolor{black}{A more detailed discussion and a kinetics interpretation of the model is provided in Section \ref{sec:discussion}.}

Consider the simple simplex
\[
	\Delta_n = \left\{ P | p_i \geq 0 , \sum \nolimits_i p_i = 1 \right\}.
\]
$\Delta_n$ is clearly forward invariant under the dynamics \eqref{eq:sysP} since it preserves non-negativity and obeys the ``conservation law'' $\sum_i p_i = \mathrm{const}$. (The latter follows immediately from the fact that $K$ has zero column sums.) Thus any solution $P(\cdot; t_0,P_0)$ of \eqref{eq:sysP} starting from $P_0=P(t_0) \in \Delta_n$ remains in $\Delta_n$ for all $t \geq t_0$.

The invariance of $\Delta_n$ under \eqref{eq:sysP} can be used to prove certain important properties of $K$ and its associated system \eqref{eq:sysP}. Two examples are presented below.
\begin{itemize}
\item {\bf Equilibria.} The non-negative vector $P^*$ such that $K P^* = 0$ is known as the Perron vector of $K$,  and defines an equilibrium of system \eqref{eq:sysP}. The existence of this vector $P^*$ can also be deduced from the forward invariance of $\Delta_n$. Indeed, as any continuous map $\Phi:\Delta_n \rightarrow \Delta_n$ has a fixed point (Brouwer fixed point theorem), $\Phi = \exp(Kt)$ has a fixed point in $\Delta_n$ for any $t\geq t_0$. If $\exp(Kt)P^* = P^*$ for some $P^* \in \Delta_n$ and sufficiently small $t>t_0$, then $K P^*=0$ because
\[
	\exp(Kt)P = P + tK P + o(t^2).
\]
\item {\bf Eigenvalues of $K$.} It is clear that $K$ has a zero eigenvalue. In fact, Gershgorin's theorem implies that all eigenvalues of $K$ are in the union of closed discs
\[
	D_i = \{ \lambda \in \mathbb{C} | \|\lambda - k_{ii} \| \leq |k_{ii}| \}.
\]
Thus $K$ does not have purely imaginary eigenvalues. This can also be deduced from the forward invariance of $\Delta_n$ in combination with the assumption of a positive equilibrium $P^*$. We exclude the eigenvector corresponding to the zero eigenvalue and consider $K$ on the invariant hyperplane where $\sum_i p_i=0$. If $K$ has a purely imaginary eigenvalue $\lambda$, then there exists a 2$D$ $K$-invariant subspace $U$, where $K$ has two conjugated imaginary eigenvalues, $\lambda$ and $\overline{\lambda}=-\lambda$. Restriction of $\exp(Kt)$ on $U$ is a one-parametric group of rotations. For the positive equilibrium $P^*$ the intersection $(U+P^*)\cap \Delta_n$ is a convex polygon. It is forward invariant with respect to \eqref{eq:sysP} because $U$ is invariant, $P^*$ is an equilibrium and $\Delta_n$ is forward invariant. But a polygon on a plane cannot be invariant with respect to the one-parametric semigroup of rotations $\exp(Kt)$ ($t\geq 0$). This contradiction proves the absence of purely imaginary eigenvalues.
\end{itemize}

The main result of this section is the following theorem:
\begin{thrm}\label{theorem:linear}
For every nonzero eigenvalue $\lambda$ of matrix $K$
\begin{equation}
\frac{|\Im \lambda |}{|\Re \lambda |} \leq
\cot\frac{\pi}{n}.
\end{equation}
\end{thrm}
The proof of this theorem can be extracted from the general Dmitriev--Dynkin--Karlelevich theorems \cite{DD1946,Karp1951}, but the straightforward geometric proof presented below, which makes use of the forward invariance of $\Delta_n$ under the dynamics \eqref{eq:sysP}, seems to be more instructive.

\begin{proof}
Let us assume that system \eqref{eq:sysP} has a positive equilibrium
$P^*\in \Delta_n$    ($p_i^*>0$ for all $i=1,\ldots,n$). For this $P^*$,
$$\sum_j q_{ij}p^*_j= \sum_j q_{ji} p^*_i .$$
Systems \eqref{eq:sysP} without strictly positive equilibria (but with non-negative ones) may be considered as limits of those with
positive equilibria.

Let $\lambda$ be a complex eigenvalue of $K$ and let $U$ be a 2D real subspace of the
hyperplane $\sum_i p_i=0$ that corresponds to the pair of complex conjugated eigenvalues,
$(\lambda, \overline{\lambda})$. Let us select a coordinate system in the plane $U+P^*$
with the origin at $P^*$ such that restriction of $K$ on this plane has the following
matrix
$$\mathcal{K}=\left[\begin{array}{cc}
\Re \lambda &- \Im \lambda\\
\Im \lambda & \Re \lambda
\end{array}\right]\, .
$$

In this coordinate system $$\exp (t\mathcal{K})=\left[\begin{array}{cc}
\exp(t \Re \lambda) \cos (t \Im \lambda) &-\exp(t \Re \lambda) \sin(t \Im \lambda)\\
\exp(t \Re \lambda) \sin(t \Im \lambda) & \exp(t \Re \lambda) \cos (t \Im \lambda)
\end{array}\right]\, .
$$

The intersection $\mathcal{A}=(U+P^*)\cap \Delta_n$ is a polygon. It has no more than
$n$ sides because $\Delta_n$ has $n$ $(n-2)$-dimensional faces (each of them is given in
$\Delta_n$ by an equation $p_i=0$). For the transversal intersections (the generic case)
this is obvious. Non-generic situations can be obtained as limits of generic cases when
the subspace $U$ tends to a non-generic position. This limit of a sequence of polygons
cannot have more than $n$ sides if the number of sides for every polygon in the sequence
does nor exceed $n$.

Let the polygon $\mathcal{A}$ have $m$ vertices $\mathbf{v}_j$ ($m\leq n$). We move the
origin to $P^*$ and enumerate these vectors  $\mathbf{x}_i=\mathbf{v}_i-P^*$
anticlockwise (Fig~\ref{Fig:Polygon}). Each pair of vectors
$\mathbf{x}_i,\mathbf{x}_{i+1}$ (and $\mathbf{x}_m,\mathbf{x}_{1}$) form a triangle with
the angles $\alpha_i$, $\beta_i$ and $\gamma_i$, where $\beta_i$ is the angle between
$\mathbf{x}_i$ and $\mathbf{x}_{i+1}$, and $\beta_m$ is the angle between $\mathbf{x}_m$
and $\mathbf{x}_{1}$. The Sine theorem gives $\frac{|\mathbf{x}_i|}{\sin
\alpha_i}=\frac{|\mathbf{x}_{i+1}|}{\sin \gamma_i}$, $\frac{|\mathbf{x}_m|}{\sin
\alpha_m}=\frac{|\mathbf{x}_{1}|}{\sin \gamma_1}$.

Several elementary identities and inequalities hold:
\begin{equation}\label{Cond}
\begin{split}
&0 <\alpha_i,\beta_i, \gamma_i<\pi; \;\; \sum_i \beta_i=2\pi; \;\; \alpha_i+ \beta_i + \gamma_i =\pi; \\
&\prod_i \sin \alpha_i =\prod_i \sin \gamma_i \mbox{ (the closeness condition).}
\end{split}
\end{equation}
These  conditions  (\ref{Cond})  are necessary and sufficient for the existence of a
polygon $\mathcal{A}$ with these angles which is star-shaped with respect to the origin.

Let us consider anticlockwise rotation ($\Im \lambda <0$, Fig.~\ref{Fig:Polygon}).
(The case of clockwise rotations differs only in notation.) For the angle $\delta$ between
$K\mathbf{x}_i$ and $\mathbf{x}_i$, $\sin \delta=-\Im \lambda$,  $\cos \delta=-\Re
\lambda$ and $\tan \delta=\frac{\Im \lambda}{\Re \lambda}$.

\begin{figure}
\centering{
\includegraphics[width=0.3\textwidth]{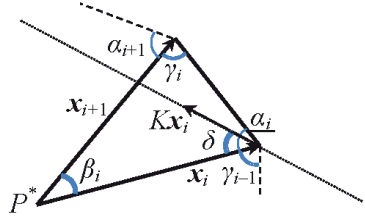}
\caption{\label{Fig:Polygon}The polygon $\mathcal{A}$ is presented as a sequence of vectors $\mathbf{x}_i$.
The angle $\beta_i$ between vectors $\mathbf{x}_i$ and $\mathbf{x}_{i+1}$ and the angles $\alpha_i$ and $\gamma_i$ of the triangle with
sides  $\mathbf{x}_i$ and $\mathbf{x}_{i+1}$ are shown.
In the Fig., rotation goes anticlockwise, i.e. $\Im \lambda <0$. In this case, the polygon $\mathcal{A}$
is invariant with respect to the semigroup $\exp (t\mathcal{K})$ ($t\geq 0$)
if and only if $\delta \leq \alpha_i$ for all $i=1, \ldots, m$, where $\delta$ is the angle between the vector field $\mathcal{K} \mathbf{x}$ and
the radius-vector $\mathbf{x}$. } }
\end{figure}

For each point $\mathbf{x} \in U+P^*$ ($\mathbf{x} \neq P^*$), the straight line
$\{\mathbf{x}+\epsilon \mathcal{K}\mathbf{x} \, | \, \epsilon \in \mathbb{R}\}$ divides
the plane $U+P^*$ in two half-plane (Fig.~\ref{Fig:Polygon}, dotted line). Direct
calculation shows that the semi-trajectory $\{\exp (t\mathcal{K})\mathbf{x}\, |\, t\geq
0\}$ belongs to the same half-plane as the origin $P^*$ does. Therefore, if $\delta \leq
\alpha_i$ for all $i=1,\ldots,m$ then the polygon $\mathcal{A}$ is forward-invariant with
respect to the semigroup $\exp (t\mathcal{K})$ ($t\geq 0$). If $\delta >\alpha_i$ for
some $i$  then for sufficiently small $t>0$  $\exp (t\mathcal{K})\mathbf{x}_i \notin
\mathcal{A}$ because $\mathcal{K}\mathbf{x}_i $ is the tangent vector to the
semi-trajectory at $t=0$. Thus, the polygon $\mathcal{A}$ is forward-invariant with respect to the semigroup $\exp (t\mathcal{K})$ ($t\geq 0$) if and only if $\delta \leq \alpha_i$ for all $i=1,\ldots,m$. The maximal $\delta$ for which $\mathcal{A}$ is still forward-invariant is $\delta_{\max}=\min_i \{\alpha_i\}$. We have to find the polygon with
$m\leq n$ and the maximal value of $\min_i \{\alpha_i\}$. Let us prove that this is a
regular polygon with $n$ sides. Let us find the maximizers  $\alpha_i,\beta_i, \gamma_i$
($i=1, \ldots, m$) for the optimization problem:
\begin{equation}\label{OptimProbl}
\min_i \{\alpha_i\} \to \max\; \mbox{subject to conditions (\ref{Cond}).}
\end{equation}
The solution of this problem is that all $\alpha_i$ are equal. To prove this equality, observe that $\min_i \{\alpha_i\}<\frac{\pi}{2}$ under conditions (\ref{Cond}) (if all
$\alpha_i\geq \frac{\pi}{2}$ then the polygonal chain $\mathcal{A}$ cannot be closed).
Let $\min_i \alpha_i=\alpha$. Let us substitute in (\ref{Cond}) the variables $\alpha_i$
which take this minimal value by $\alpha$. The derivative of the left hand part of the
last condition in (\ref{Cond}) with respect to $\alpha$ is not zero because
$\alpha<\frac{\pi}{2}$. Assume that there are some $\alpha_j> \alpha$. Let us fix the
values of $\beta_i$ ($i=1, \ldots, m$). Then $\gamma_i$ is a function of $\alpha_i$,
$\gamma_i=\pi-\beta_i-\alpha_i$. We can use the implicit function theorem to increase
$\alpha$ by a sufficiently small number $\varepsilon>0$ and to change the non-minimal
$\alpha_j$ by a small number too, $\alpha_j\mapsto \alpha_j-\theta$;
$\theta=\theta(\varepsilon)$. Therefore, at the solution of (\ref{OptimProbl}) all
$\alpha_j = \alpha$ ($j=1,\ldots,m$).

Now, let us prove that for solution of the problem (\ref{OptimProbl}) all $\beta_i$ are
equal. We exclude $\gamma_i$ from conditions (\ref{Cond}) and write $\beta_i+\alpha<
\pi$; $0<\beta_i,\alpha$;
\begin{equation}\label{NewCOnd}
m \log\sin \alpha =\sum_i \log \sin (\beta_i+\alpha).
\end{equation}
Let us consider this equality as equation with respect to unknown $\alpha$. The function
$\log \sin x$ is strictly concave on $(0,\pi)$. Therefore, for $x_i\in (0,\pi)$
$$\log \sin\left( \frac{1}{m}\sum_{i=1}^m x_i\right)\geq \frac{1}{m}\sum_{i=1}^m\log \sin x_i$$
and the equality here is possible only if all $x_i$ are equal.  Let $\alpha^*\in
(0,\pi/2)$ be a solution of (\ref{NewCOnd}). If not all the values of $\beta_i$ are equal
and we replace $\beta_i$ in (\ref{NewCOnd}) by the average value, $\beta=\frac{2\pi}{m}$,
then the value of the right hand part of (\ref{NewCOnd}) increases and $\sin \alpha^* <
\sin (\beta+\alpha^*)$. If we take all the $\beta_i$ equal then (\ref{NewCOnd})
transforms into elementary trigonometric equation $\sin \alpha = \sin (\beta+\alpha)$.
The solution $\alpha$ of equation (\ref{NewCOnd}) increases when we replace $\beta_i$  by
the average value: $\alpha> \alpha^*$ because $\sin \alpha^* < \sin (\beta+\alpha^*)$,
$\alpha \in (0,\pi/2)$ and $\sin \alpha$ monotonically increases on this interval.  So,
for the maximizers of the conditional optimization problem (\ref{OptimProbl}) all
$\beta_i=\frac{2\pi}{m}$ and $\alpha_i=\gamma_i=\frac{\pi}{2}-\frac{\pi}{m}$. The maximum
of $\alpha$ corresponds to the maximum of $m$. Therefore, $m=n$. Finally, $\max
\{\delta\}=\frac{\pi}{2}-\frac{\pi}{n}$ and
$$\max\left\{\frac{|\Im \lambda |}{|\Re \lambda |}\right\} =
\cot\frac{\pi}{n} .$$
\end{proof}

\begin{rmrk} It is important to note that the bound given in Theorem \ref{theorem:linear} is sharp. Indeed, let $K$ define a directed cycle with uniform weights $q$, e.g.
\[
	K = \left(\begin{array}{rrrrrrr}
	-q & 0 & 0 & \cdots & 0 & q \\
	q & -q & 0 & \cdots & 0 & 0 \\
	0 & \ddots & \ddots & \ddots & \vdots & \vdots \\
	\vdots & \ddots & \ddots & \ddots & 0 & 0 \\
	0 &  & \ddots & q & -q & 0 \\
	0 & 0 & \cdots & 0 & q & -q \\
	\end{array}\right).
\]
The eigenvalues of $K$ are
\[
	\lambda_k = -q + q \exp\left(\frac{2 \pi k \mathrm{i} }{n} \right),\quad k=0,1,\ldots,n-1,
\]
cf. \cite{davis}, with $\mathrm{i}=\sqrt{-1}$ the imaginary unit. Thus $\frac {|\mathfrak{J} \lambda_1|}{ |\mathfrak{R} \lambda_1|} = \cot\frac{\pi}{n}$. Note that for large $n$,
\[
\cot\frac{\pi}{n} \approx \frac{n}{\pi},
\]
which means that oscillations in a simple cycle with a large number of systems decay very slowly.
\end{rmrk}

An important consequence of this extremal property of a simple cycle is that not only transients in the cycle decay very slowly but also that the overall behavior of transients becomes extremely sensitive to perturbations. This, as we show in the next sections, gives rise to resonances and bistabilities if neutrally stable nodes in (\ref{eq:sysP}) are replaced with ones exhibiting oscillatory dynamics. As a model of nodes with oscillatory activity the classical Fitzhugh-Nagumo \cite{FitzHugh} system has been chosen. Our choice of this system among various alternatives \cite{Izhikevich} was motivated purely by its simplicity and relevance for modelling behavior of neural systems.

\setcounter{equation}{0}

\section{Coupled Nonlinear Neural Oscillators}\label{sec:neural}

Consider a network of FitzHugh-Nagumo (FHN) neurons
\begin{equation}\label{eq:FHN}
	\left\{ \begin{array}{l}
		\dot z_j = \alpha \left( y_j - \beta z_j \right) \\
		\dot y_j = y_j - \gamma y_j^3 -z_j + u_j,
	\end{array}\right.
\end{equation}
$j=1,2,\ldots,n$ with parameters $\alpha,\beta,\gamma$ chosen as
\[
\alpha=\tfrac{8}{100}, \ \beta=\tfrac{8}{10}, \ \gamma=\tfrac{1}{3}
\]
The FHN neurons interact via diffusive coupling
\begin{equation}\label{eq:coupling}
	u_j = \sigma \sum_{ l=1}^n q_{j l} (y_l-y_j)
\end{equation}
with constant $\sigma \in\mathbb{R}$, $\sigma>0$, being the coupling strength.
For convenience, let
\[
y=(y_1,\dots,y_n), \ u=(u_1,\dots,u_n),\ x=(y,z),
\]
and $x(\cdot;x_0,\sigma)$ denote a solution of the coupled system with the coupling strength $\sigma$ and satisfying the initial condition $x(0)=x_0$. The topology of network connections in (\ref{eq:coupling}) is characterized by the adjacency matrix $Q$ with zeros on the main diagonal and entries identical to the values of $q_{ij}$ for $i\neq j$, $i,j\in\{1,\dots,n\}$. The matrix $Q$ is now assumed to be a circulant matrix
\[
Q=\left(\begin{array}{ccccc} 0 & 0 & \cdots & 0 &1\\
                             1 & 0 & \cdots & 0 & 0\\
                             0 & 1 & \ddots & \vdots & \vdots \\
                             \vdots & \ddots & \ddots & 0 & 0\\
                             0 & \cdots & 0 & 1  & 0\end{array}\right).
\]
Thus besides assuming the network structure to be a directed ring we have also assumed the interaction weights $q_{jl}$ to be identical and, without loss of generality, we have set these weights of interaction to $1$.

At first glance, the connectivity pattern specified by $Q$ differs from that specified by matrix $K$ in (\ref{eq:K}). Yet, if coupling (\ref{eq:coupling}) is rewritten in the vector-matrix notation then the following identity holds
\begin{equation}\label{eq:coupling:vector}
u =  \sigma (Q-I_n) y \triangleq - \sigma L y.
\end{equation}
As remarked before, the network Laplacian matrix
\[
L= \mathrm{Diag}(\sum_{j\neq l}q_{jl})-Q=I_n-Q
\]
can be related to the matrix $K$ corresponding to the simple cycle in ``reverse'' direction as $L=-K^T$.

In what follows we will employ the notions of {\it semi-passivity} and {\it strict semi-passivity} that have been introduced first in \cite{Pogromsky1998}. For consistency, we recall these notions below

\begin{dfntn} \it Consider a system of first-order nonlinear ordinary differential equations
\[
\dot{x}=f(t,x,u(t)),
y=h(x)
\]
where $f:\Real\times\Real^n\times \Real\rightarrow\Real^n$ is a continuous and locally Lipschitz function, $h:\Real^n\rightarrow\Real$ is a continuous function, and $u:\Real\rightarrow\Real$ is a continuous function. Let $x(\cdot;t_0,x_0,[u])$ be a solution of the Cauchy problem $x(t_0;t_0,x_0,[u])=x_0$, and let $\mathcal{U}\subset\mathcal{C}^0$ be the set of inputs $u$ for which the solution $x(\cdot;t_0,x_0,[u])$ is defined in forward time.

The system is called semi-passive if there is a non-negative function $S:\Real^n\rightarrow\Real_{+}$  (a storage function) and  a function $H:\Real^{n}\rightarrow \Real$ such that for each $x(\cdot;t_0,x_0,[u])$ the following holds for all $t\geq t_0$ in the domain of this solution definition:
\[
S(x(t;t_0,x_0,[u]))-S(x_0)\leq \int_{t_0}^{t} y(\tau)u(\tau) - H(x(\tau,t_0,x_0,[u])) d\tau,
\]
where the function $H$ is non-negative outside a ball in $\Real^n$.

The system is called strictly semi-passive if the function $H$ is strictly positive outside a ball in $\Real^n$.
\end{dfntn}

\subsection{Boundedness of solutions in the coupled system}

\begin{lmm}\label{lem:boundedness} The solutions of the ring network of FHN neurons are ultimately bounded uniformly in $x_0$, $\sigma\in\Real_{\geq 0}$. That is, there is a compact set $\Omega\in\Real^{2n}$ such that for all $x_0\in\Real^{2n}$, $\sigma\in\Real_{\geq 0}$
\[
\lim_{t\rightarrow\infty} \mathrm{dist}\left(x(t,x_0,\sigma),\Omega\right)=0.
\]
\end{lmm}

\begin{proof} We being begin with establishing that the FHN neuron is strictly semi-passive (see also \cite{Steur2009}).

Let $S(z_j,y_j) = \tfrac{1}{2} \left( \alpha^{-1} z_j^2 + y_j^2 \right)$ be the storage function. Then
\[
	\dot S = -H(z_j, y_j) + y_j u_j
\]
with $H(z_j,y_j) = \beta z_j^2 + y_j^2 \left( \gamma y_j^2-1 \right)$. Noticing that
\begin{equation}\label{eq:semi-passivity:dV}
\begin{split}
&\beta z_j^2 + y_j^2 \left( \gamma y_j^2-1 \right)=\beta z^2 + d y^2 + \gamma y^4 - y^2 - dy^2=\\
&\beta z^2 + d y^2 + \left(\sqrt{\gamma}y^2-\tfrac{d+1}{2 \sqrt{\gamma}}\right)^2 -\tfrac{(d+1)^2}{4\gamma}
\end{split}
\end{equation}
we can conclude that $H(z_j,y_j)$ is positive for all $z_j,y_j$ such that
\begin{equation}\label{eq:semi-passivity:ellipse}
\beta z_j^2 + d y_j^2 > \tfrac{(d+1)^2}{4\gamma}.
\end{equation}
Assigning the value of $d$ in (\ref{eq:semi-passivity:ellipse}) as $d=\beta$, ensures that $H(z_j,y_j)$ is positive outside the ball
 \[
 z_j^2+y_j^2\leq \tfrac{(\beta+1)^2}{4\beta \gamma}.
 \]

Now consider $V(z,y) = S(z_1,y_1) + \ldots + S(z_n, y_n)$. Then the strict semi-passivity property of the FHN neurons implies
\[
	\dot V \leq -H(z_1,y_1) - \ldots - H(z_n,y_n) - \sigma y^T L y.
\]
Notice that the matrix $L+L^T$ is the Laplacian matrix of the undirected ring, which is known to be positive semi-definite. Hence
\[
	y^T L y = \tfrac{1}{2} y^T (L+L^T) y \geq 0,
\]
and consequently\\
\[
	\dot V \leq -H(z_1,y_1) - \ldots - H(z_n,y_n).
\]
Therefore, setting the value of $d$ in (\ref{eq:semi-passivity:dV}) equal to $\alpha\beta$ results in
\[
\dot V \leq - \sum_{j=1}^n \beta z_j^2 + \beta \alpha y_j^2 + n \tfrac{(\alpha\beta
+1)^2}{4\gamma} = - \beta\alpha V + n \tfrac{(\alpha\beta
+1)^2}{4\gamma}.
\]
Noticing that the function $V$ is radially unbounded, positive-definite, we invoke the Comparison Lemma (see e.g. \cite{Khalil}) in order to conclude that solutions of the coupled system are bounded and converge asymptotically to a compact set of which the size is independent of the parameter $\sigma$.
\end{proof}

\subsection{Sufficient conditions for synchronization}

\subsubsection{Directed chain: "no looking back"}
First we consider the dynamics of two coupled systems in the leader-follower configuration:
 \begin{align}\label{eq:directed_chain}
	&\left\{ \begin{array}{l}
		\dot z_1 = \alpha \left( y_1 - \beta z_1 \right) \\
		\dot y_1 = y_1 - \gamma y_1^3 -z_1
	\end{array}\right.  \\[2ex]
	&\left\{ \begin{array}{l}
		\dot z_2 = \alpha \left( y_2 - \beta z_2 \right) \\
		\dot y_2 = y_2 - \gamma y_2^3 -z_2 +\sigma(y_1 - y_2).
	\end{array}\right. \nonumber
 \end{align}
\begin{thrm}\label{theorem:chain} Consider the system of coupled FHN oscillators (\ref{eq:directed_chain}) in which the parameter $\sigma$ is chosen so that
\[
\sigma>1.
\]
Then solutions of the system asymptotically synchronize for all values of initial conditions.
\end{thrm}

 \begin{proof}
In accordance with Lemma \ref{lem:boundedness} solutions of the coupled system exist and are bounded for all $t>0$. Define
\[
	\tilde z = z_1 - z_2, \quad \tilde y = y_1-y_2,
\]
such that
\begin{align*}
	&\dot{\tilde z} = \alpha\left(\tilde y - \beta \tilde z \right) \\
	&\dot{\tilde y} = \tilde y -\gamma (y_1^3 - y_2^3) - \tilde z - \sigma \tilde y.
\end{align*}
Consider the function
\[
	V = \tfrac{1}{2} \left(\tfrac{1}{\alpha} \tilde z^2 + \tilde y^2 \right),
\]
then, using the equality
\[
(y_1-y_2)(y_1^3 - y_2^3) = \tfrac{1}{4} (y_1 -y_2)^2 \left( 3 (y_1+y_2)^2 + (y_1 -y_2)^2 \right),
\]
we find
\[
	\dot V = -\beta \tilde z^2  + (1-\sigma) \tilde y^2 -\tfrac{\gamma}{4} \tilde y^2 \left(3 (y_1 +y_2)^2+ \tilde y^2) \right).
\]
Thus if $\sigma>1$ we have $\dot V<0$ and the chain of FHN neurons synchronizes.
\end{proof}

Generalizing two coupled systems to a directed chain of $n$ oscillators, we  observe that the Laplacian matrix of this configuration is
 \[
 	L = \left( \begin{array}{rrrrr}
 	0 & 0 & 0& \cdots & 0 \\
 	-1 & 1 & 0& \cdots & 0 \\
 	0 & -1 & 1 & \ddots & \vdots \\
 	\vdots & \ddots & \ddots & \ddots & 0 \\
 	0 & \cdots & 0 & -1 & 1 \\
 	\end{array}\right).
 \]
The matrix $L$ has only real eigenvalues; a simple zero eigenvalue and $n-1$ eigenvalues equal to $1$. The only type of stable correlated oscillations we can find in the chain are the completely synchronous oscillations. These synchronous oscillations will emerge for values of the coupling strength $\sigma$ for which the chain of $2$ FHN oscillators synchronize. Thus the conditions for synchronization are independent of the size of the network (i.e. the length of the chain).  Numerical simulations below illustrate this statement.

\begin{figure}
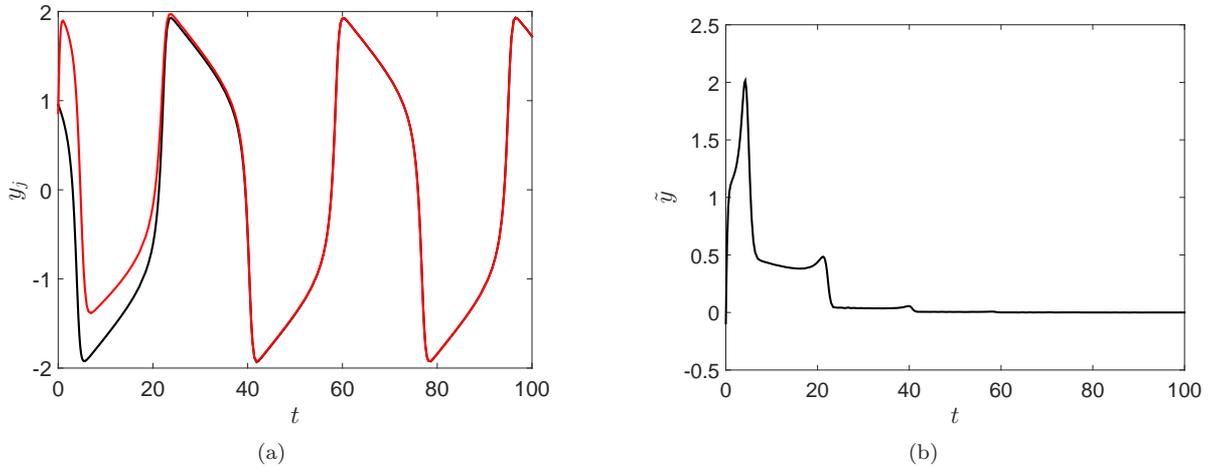

\subfigure[]{
\includegraphics[width=0.45\textwidth]{fhnchain2.eps}
}\hfill
\subfigure[]{
\includegraphics[width=0.45\textwidth]{fhnchain2errors.eps}
}
\caption{Synchronization of two FHN oscillators for $\sigma=1.5$. (a) Outputs of FHN oscillator 1 (leader, black) and FHN oscillator 2 (follower, red). (b) Synchronization output error $\tilde y:= y_1-y_{2}$.}\label{fig:LF1}
\end{figure}

\begin{figure}
\subfigure[$n=10$]{
\includegraphics[width=0.45\textwidth]{fhnchain10errors.eps}
}\hfill
\subfigure[$n=50$]{
\includegraphics[width=0.45\textwidth]{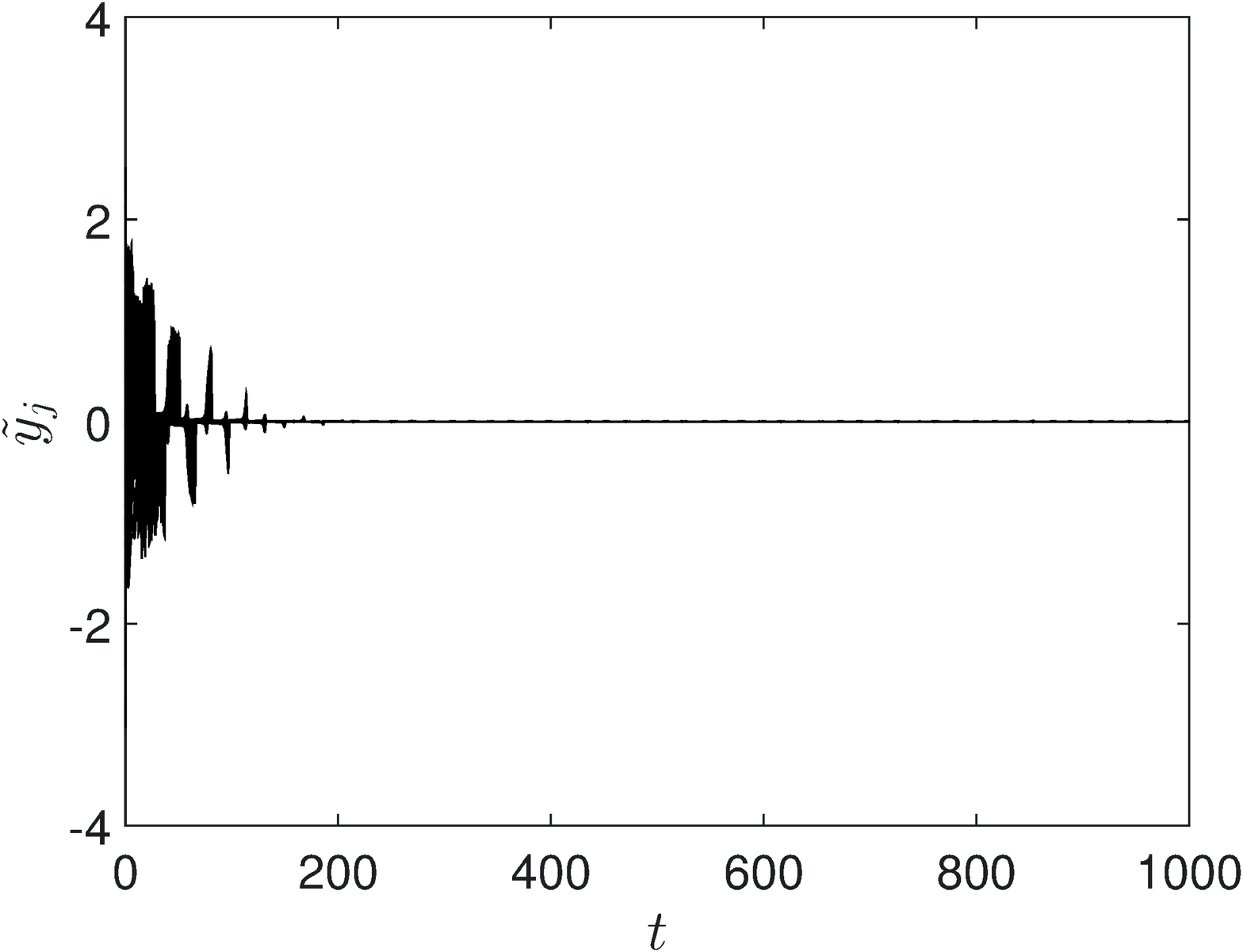}
}
\subfigure[$n=100$]{
\includegraphics[width=0.45\textwidth]{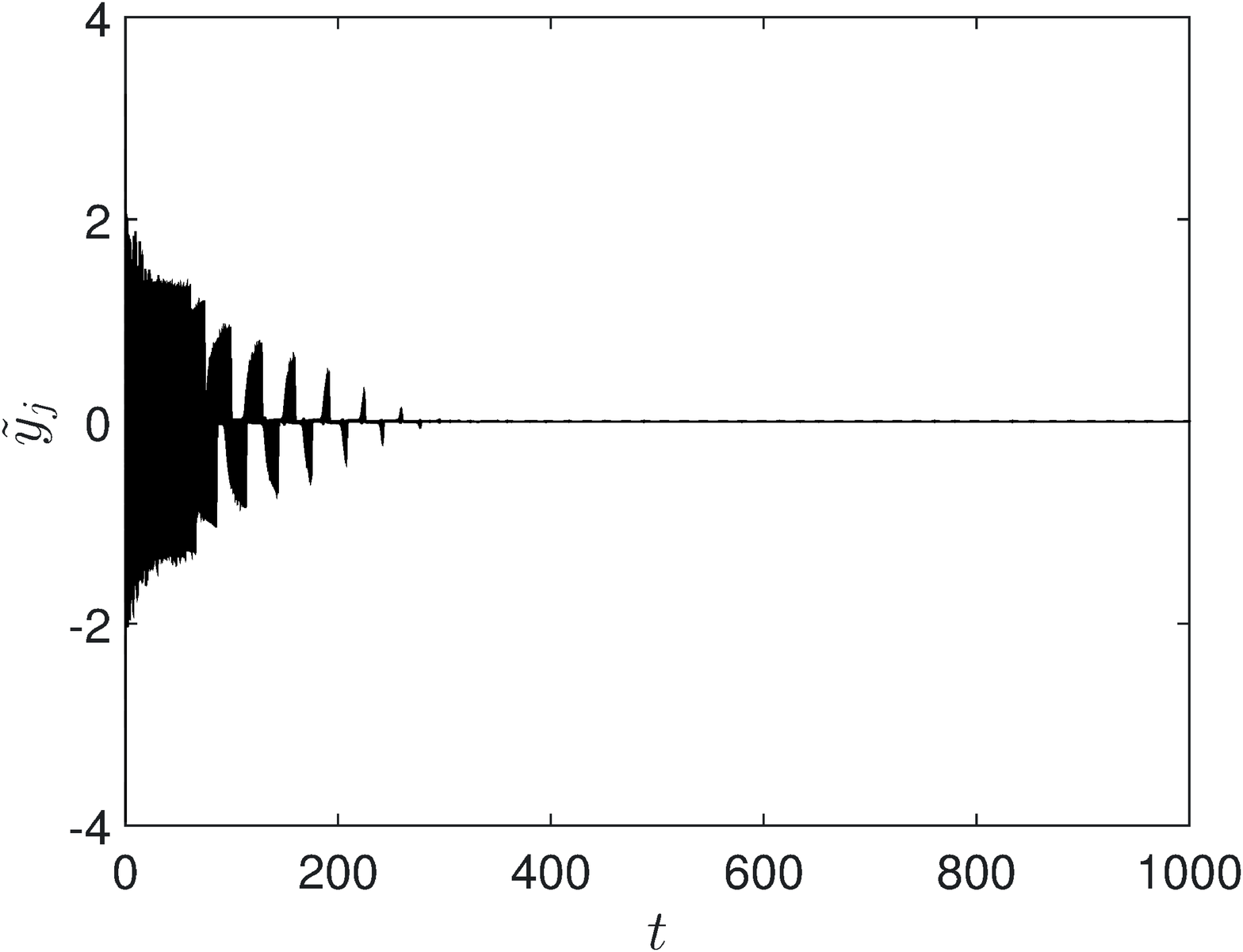}
}\hfill
\subfigure[$n=150$]{
\includegraphics[width=0.45\textwidth]{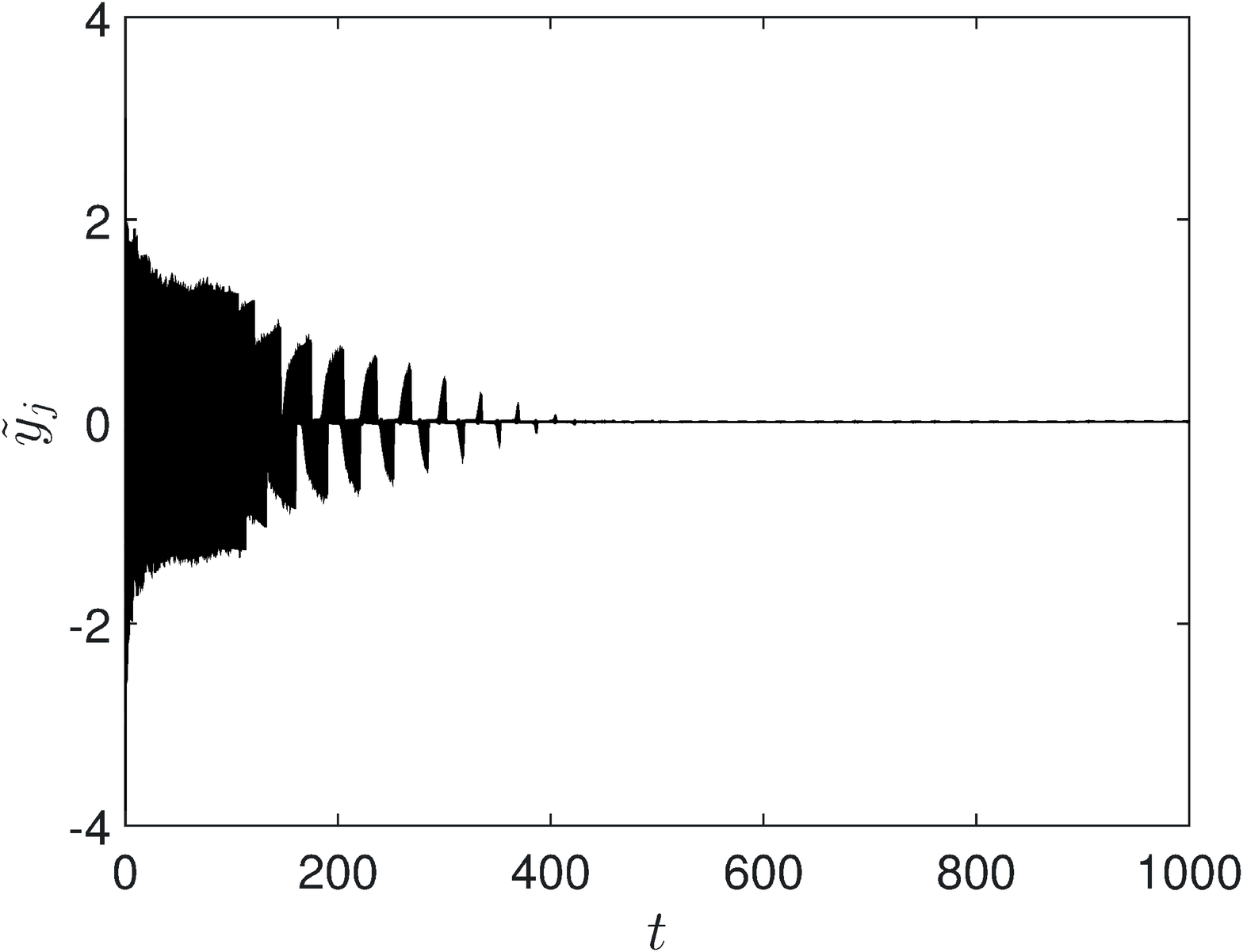}
}
\caption{Synchronization output errors $\tilde y_j:= y_j-y_{j-2}$, $j=2,\ldots,n$, for $\sigma=1.5$ and different lengths of the chain.}\label{fig:LF2}
\end{figure}

Figure \ref{fig:LF1} shows the outputs of two FHN oscillators and the synchronization output error for $\sigma=1.5$. Figure \ref{fig:LF2} shows the results for longer chains; Even though the convergence to the synchronous state is slower for longer chains, the oscillators in the chains always end up in synchrony.

\subsubsection{Directed ring: "looking back"}

Suppose now that the $n$-th oscillator is feeding back its output to the input of the $1$st, that is the network topology is that of the directed ring. As we shall see later the presence of such an extra connection has a drastic effect on the system's performance with respect to the coupling strength needed to maintain stable full-state synchrony. This is reflected in the statement of the theorem below.

\begin{thrm}\label{theorem:ring}
Consider the system of coupled FHN oscillators (\ref{eq:FHN}), (\ref{eq:coupling}), and let $\lambda_1\leq\lambda_2\leq\dots\leq\lambda_n$ be the eigenvalues of the symmetrized Laplacian of the network $\tfrac{1}{2}(L+L^T)$. Then solutions of the coupled system asymptotically synchronize providing that
\[
\sigma  \lambda_2 > 1
\]
\end{thrm}
\begin{proof}
Consider the new variables
\[
	\tilde z = L z,\quad  \tilde y = L y,
\]
where $L$ the Laplacian matrix of the ring, i.e.
\[
\tilde z = \begin{pmatrix}
\tilde z_1 \\ \tilde z_2 \\ \vdots \\ \tilde z_n
\end{pmatrix}
=\begin{pmatrix}
z_{1}-z_{n} \\ z_{2}-z_{1} \\ \vdots \\ z_{n-1}-z_{n}
\end{pmatrix}
\text{~~and~~}
\tilde y = \begin{pmatrix}
\tilde y_1 \\ \tilde y_2 \\ \vdots \\ \tilde y_n
\end{pmatrix}
=\begin{pmatrix}
y_{1}-y_{n} \\ y_{2}-y_{1} \\ \vdots \\ y_{n-1}-y_{n}
\end{pmatrix}.
\]
It is clear that the systems are synchronized if and only if
\[
	\tilde z=0 \text{~~and~~}  \tilde y=0.
\]
Observe that $\boldsymbol{1} \notin \mathrm{range}(L)$, hence there exist no vectors $z$ and $y$ such that
\[
	L z = \boldsymbol{1}\text{~~and~~} L y = \boldsymbol{1}.
\]
This means that the projections of $(z,y)$ via $L$ take values in the set
\[
	\Omega := \{ (\tilde z,\tilde y) \in \mathbb{R}^{2n} | \tilde z \perp \boldsymbol{1}, \tilde y \perp \boldsymbol{1} \}.
\]
Thus all synchronization errors $\tilde z$ and $\tilde y$ are orthogonal to $\boldsymbol{1}$.

Consider the function $V:\Omega \rightarrow \mathbb{R}_+$:
\[
	V = \tfrac{1}{2} \left(\tfrac{1}{\alpha} \tilde z^T \tilde z + \tilde y^T \tilde y \right).
\]
From the discussion on synchronization in the chain it follows that
\[
	\dot V \leq -\beta \tilde z^T \tilde z + \tilde y^T (I-\sigma L) \tilde y - \tilde y^T W \tilde y
\]
where
\[
	W = \frac{\gamma}{4} \begin{pmatrix}
	 \left(3 (y_1 +y_n)^2+\tilde y_1^2 \right) & & \\ & \ddots & \\ & & \left(3 (y_{n+1} +y_n)^2+\tilde y_n^2 \right)
	\end{pmatrix},
\]
which is positive semi-definite, hence
\[
	\dot V \leq -\beta \tilde z^T \tilde z + \tilde y^T (I-\sigma L) \tilde y
\]
For all vectors $\tilde y \perp \boldsymbol{1}$ the following inequality holds true:
\[	
	\tilde y^T (\sigma L -I) \tilde y = \tilde y^T (\sigma \tfrac{1}{2}(L+L^T) -I) \tilde y \geq (\sigma \lambda_2-1) \tilde y^T \tilde y
\]
where $\lambda_2=\lambda_2(\tfrac{1}{2}(L+L^T))$ is the smallest non-zero eigenvalue of $\tfrac{1}{2}(L+L^T)$. An application of LaSalle's invariance principle, cf. \cite{LaSalle}, implies that the synchronization errors $\tilde z$ and $\tilde y$ converge to zero asymptotically.
\end{proof}

\begin{crllr}\label{corollary:nonlinear} For the network of $n$ coupled FHN oscillators, solutions globally asymptotically synchronize if the following inequality holds:
\[
\sigma \left(1-\cos\left(\frac{2\pi}{n}\right)\right)>1
\]
\end{crllr}
\begin{proof}
Note that $\tfrac{1}{2}(L+L^T)$ is the Laplacian matrix of the undirected ring, which has a simple zero eigenvalue with corresponding eigenvector in $\mathrm{span}(\boldsymbol{1})$. According to the properties of the spectrum of circulant matrices, cf. \cite{davis}, we know that the second smallest eigenvalue $\lambda_2$ of the symmetrized Laplacian $\tfrac{1}{2}(L+L^T)$ equals the real part of the smallest (in absolute value) non-zero eigenvalue of $L$, which we denote as $\Re(\lambda_2(L))$. Then if
\[
	\sigma \lambda_2(\tfrac{1}{2}(L+L^T))=\sigma \Re (\lambda_2(L)) > 1,
\]
we have $\dot V <0$, i.e. $V$ is a Lyapunov function on $\Omega$. Note that
\[
\lambda_2(L)=1-e^{\frac{2\pi i (n-1)}{n}}
\]
from which the result immediately follows.
\end{proof}

\subsection{Synchronization and rotating waves}

The results in the previous sections show that, on the one hand, when a system has a directed ring topology and the number of systems in the ring grows then their relative dynamics becomes more and more underdamped (Theorem \ref{theorem:linear}). On the other hand, in accordance with Corollary \ref{corollary:nonlinear}, estimates of attraction rates of the diagonal synchronization manifold rapidly diminish to zero with increasing numbers of systems. The latter result is, however, sufficient and may be conservative. To get a clearer view of the network dynamics we performed an exhaustive numerical exploration of the system dynamics for various values of coupling strengths $\sigma$ as well as network sizes $n$.

We construct a grid $(n,\sigma)$ for number of systems $n=2,\dots,20$ and coupling strengths $\sigma = \{0.05, 0.1, 0.15,\ldots,10\}$. For each $(n,\sigma)$ 100 sets of initial conditions are drawn uniformly randomly from the domain $\vert y_i(0)\vert \leq  \frac{3}{2}\sqrt{3}$ and $\vert z_i(0)\vert \leq \frac{15}{8}\sqrt{3}$, which can be shown to be positively invariant for both connectivity configurations (i.e. the directed simple cycle and  the directed chain). The \textit{MATLAB} numerical solver \textit{ode45} was used with relative and absolute error tolerances of order $10^{-5}$ to integrate dynamics for a maximum of $20,000$ time steps. At regular intervals of $1000,2000,\ldots,20,000$ time steps we interrupt integration to check for synchronization or rotating wave solutions. After $20,000$ time steps, if neither synchronization nor a rotating wave solution is detected, we register `no solution'.

Synchronization is identified in terms of the absolute error between the states of neighbouring systems averaged over a $1000$ time step window being less than $2\times 10^{-5}$. In case of no synchronization, we investigate the existence of rotating waves of Mode Type $1$. Rotating waves are defined as periodic solutions where all systems take identical orbits with constant non-zero and equal phase shifts between neighbouring systems. The mode type describes the group velocity of the wave; for a periodic wave, Mode Type $1$ describes  the case where the period of a rotating wave having non-zero wave velocity equals the period of individual oscillators. Identical orbits are identified if the absolute difference between the time shifted orbits - so that orbits are in-phase - of neighbouring systems averaged over the period of the orbit is less than $10^{-4}$. Constant and equal phase shifts (for a Mode Type 1 rotating wave) are identified if the maximum from all absolute differences between $n$ times the phase shifts between pairwise neighbouring systems and period $T$ is less than a tolerance of $10^{-2}$.

\begin{figure}[!h]
\centering
\includegraphics[width=260pt]{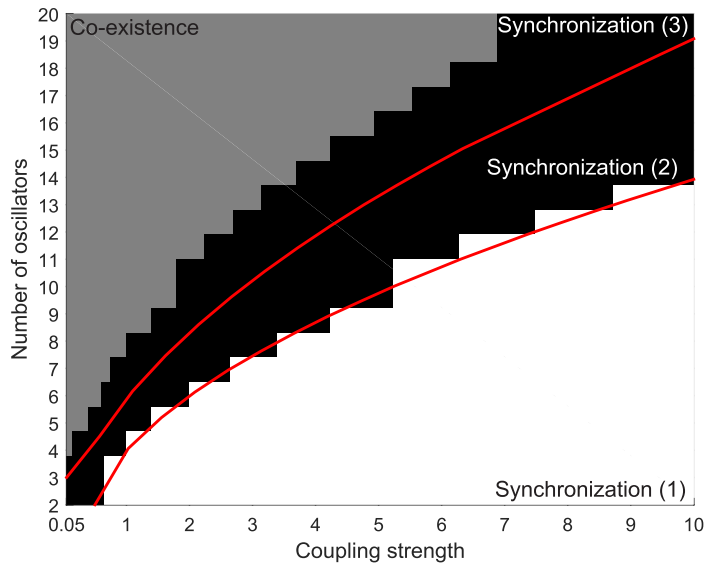}
\caption{Bifurcation diagram for directed rings of FHN oscillators. The diagram is divided into the four regions of parameter space corresponding to: \textit{Synchronization (1)}, global asymptotic synchronization that is guaranteed by the semi-passivity argument; \textit{Synchronization (2)}, synchronization registered for every set of random initial conditions during numerical simulations; \textit{Synchronization (3)}  synchronization registered for every set of random initial conditions during numerical simulations, but Floquet stability analysis of solutions of the auxiliary system indicated existence of a locally asymptotically stable rotating wave solution;  \textit{Co-existence}, both the fully synchronous and rotating wave solutions were registered during numerical simulation.}\label{Fig:FHN_network}
\end{figure}

The results of this exploration are summarized in Figure \ref{Fig:FHN_network}. This figure shows that in addition to regions corresponding to mere full asymptotic synchronization there is a wide range of parameter combinations (growing with system size) for which the system admits an asymptotically stable rotating wave solution. The larger the number of systems, the larger values of the coupling parameter $\sigma$ are required to maintain global stability of the fully synchronous state.


Two solid curves approximate boundaries between the parameter domains corresponding to analytically determined globally asymptotically stable full-state synchrony, and a partition of numerically determined globally asymptotically stable full-sate synchrony; The first region corresponds to synchronization registered for every set of random initial conditions during numerical simulations, whilst the second region corresponds to, in addition to synchronization registered for every set of random initial conditions during numerical simulations, where Floquet stability
analysis of solutions of the auxiliary system indicated existence of a locally asymptotically stable rotating wave solution.

The first (lower) curve  - separating analytical and numerical synchronization - was determined previously in the semi-passivity argument. The second (upper) curve - partitioning numerically determined globally asymptomatic synchronization  - is determined from a local stability analysis (using Floquet theory) of the rotating wave solution. Details of the second are provided below.

\begin{figure}
\centering
\includegraphics[width=260pt]{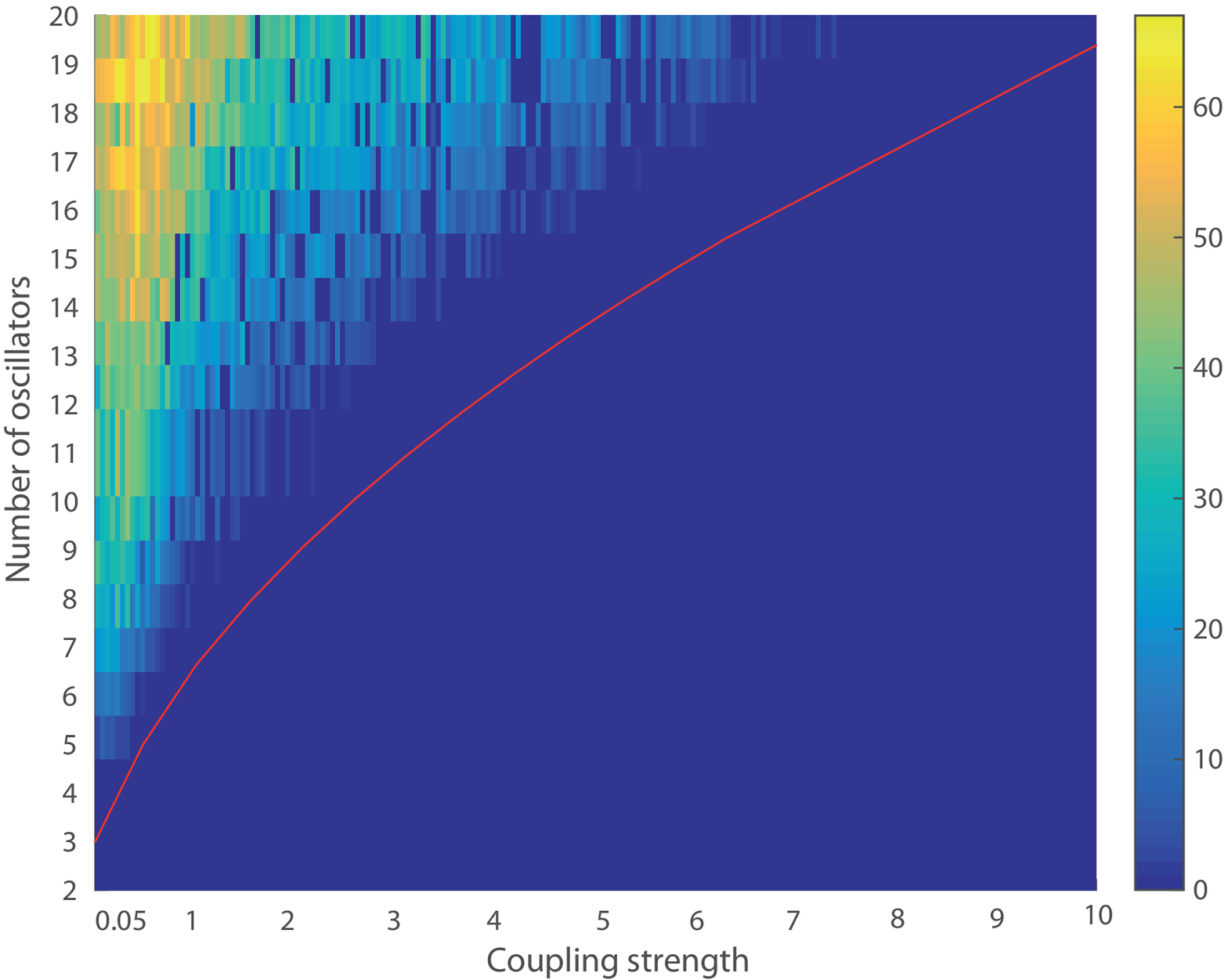}
\caption{Proportion of samples that yield a rotating wave solution of Mode Type $1$. Red curve bounds the upper left region for which Floquet stability analysis indicated the existence of a rotating wave solution.}\label{Fig:FHN_network_wave_basin_type1}
\end{figure}

\begin{figure}
\centering
\includegraphics[width=260pt]{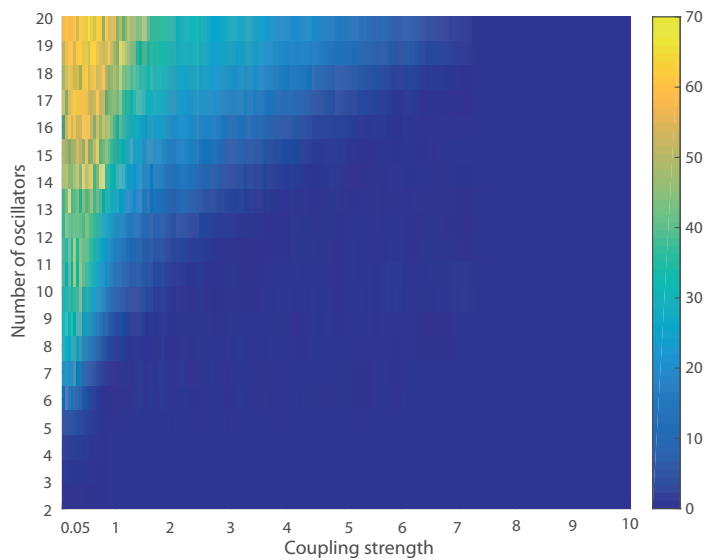}
\caption{Proportion of samples that yield a rotating wave solution for all mode types.}\label{Fig:FHN_network_wave_basin_all}
\end{figure}

Figure \ref{Fig:FHN_network_wave_basin_type1} shows for each $\sigma$ and $n$ the proportion of initial conditions that yield a rotating wave solution of Mode Type $1$ whilst Figure \ref{Fig:FHN_network_wave_basin_all} shows for all mode types, i.e. rotating waves that resonate with individual systems period of oscillation. For low coupling $\sigma$ and for increasing number of systems $n$, rotating wave solutions are found more often. This suggests a larger basin of attraction for the rotating wave than that for synchronization, and that this basin grows with increasing $n$ and decreasing $\sigma$ whilst at the same time the basin of attraction for synchronization shrinks. The relative sizes of basin of attraction result in higher or lower likelihoods for the systems to converge to a certain solution given uniformly random initial conditions.

\subsubsection{Local stability analysis of the rotating wave}
Throughout this section we consider only the rotating wave of Mode Type $1$. Similar analysis can also be performed for other mode types.

\begin{figure}
\subfigure[]{
\includegraphics[width=0.45\textwidth]{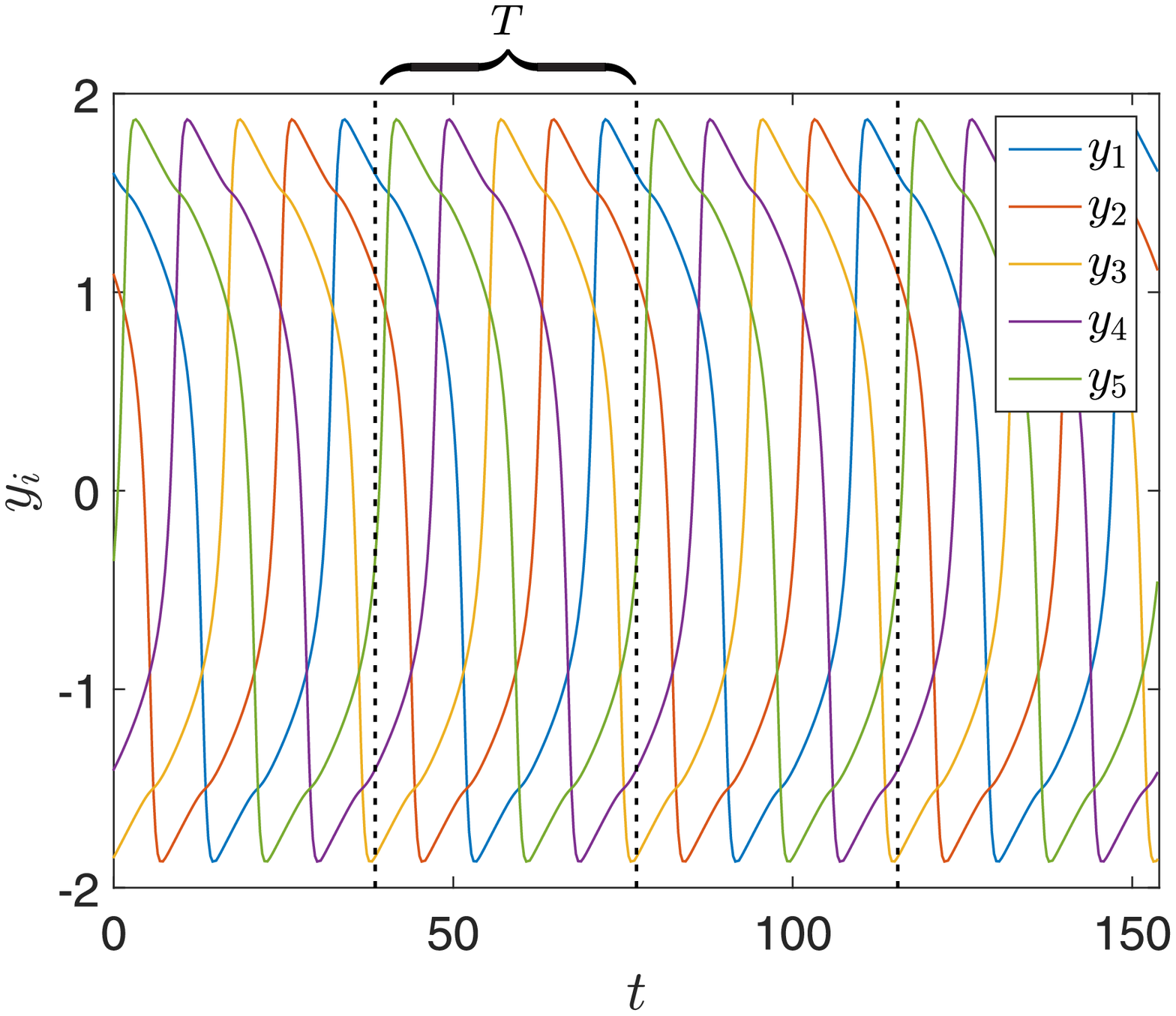}
}\hfill
\subfigure[]{
\includegraphics[width=0.45\textwidth]{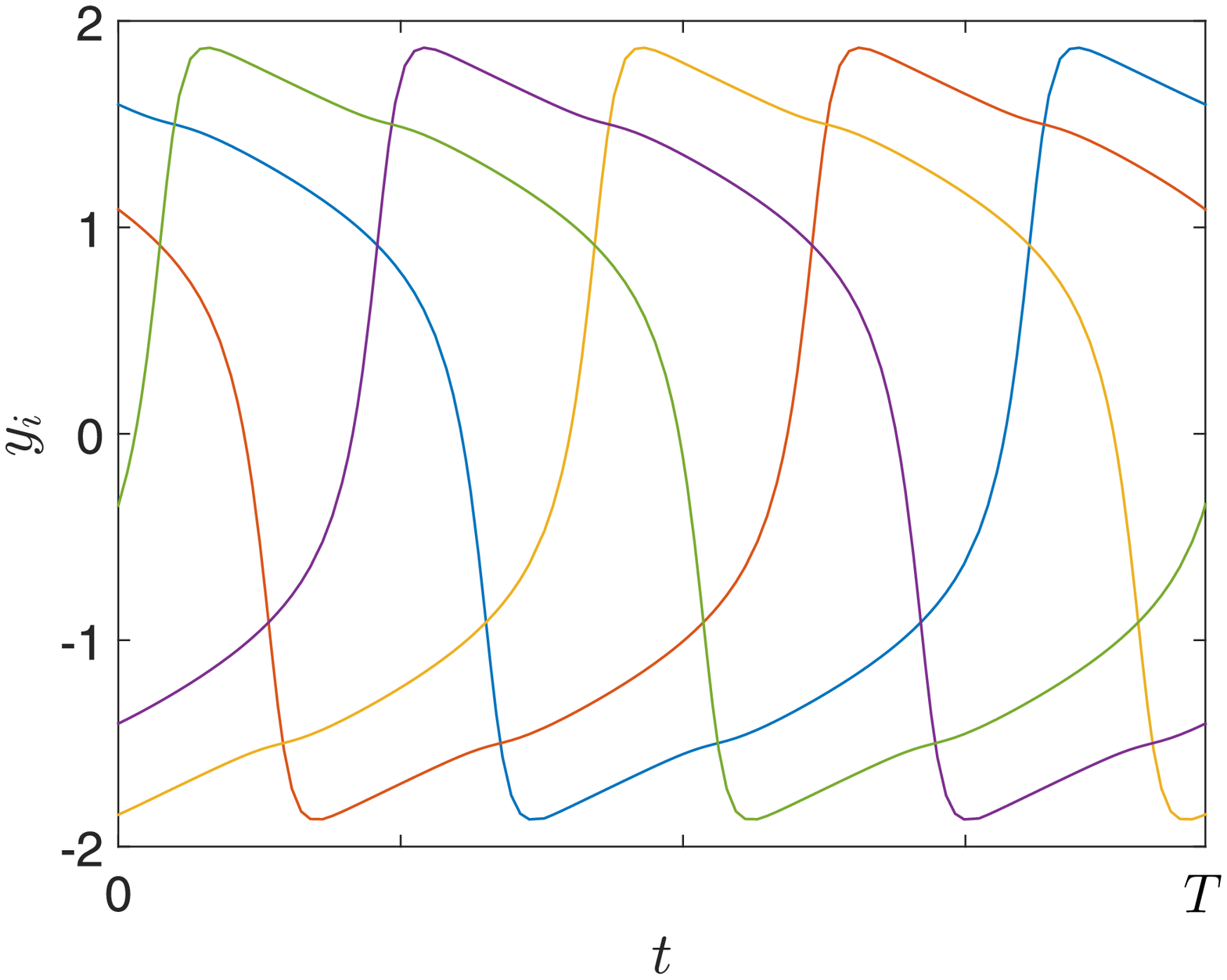}
}
\caption{(a) The $y$ dynamics over one oscillation. (b) The periodic $y$ dynamics in the time interval $[0,T]$.
The $z$ dynamics show the same type of time shifted and periodic behavior as the $y$ dynamics.}\label{fig:FHN_wave_ex}
\end{figure}

Suppose that $n$ identical coupled systems have a non-constant, $T$-periodic solution $x_j=(z_j,y_j)$ for constant $T>0$, and for which the orbit of each system is identical and time shifted by some constant $\tau = \frac{T}{n}$:
\begin{eqnarray}\label{eq:rotating_wave_solution}
x_1(t)  &=& x_2(t+\tau) = x_3(t+2\tau) = \cdots = x_n(t+(n-1)\tau)\nonumber\\
& =&  x_1(t+n\tau) = x_1(t+T).
\end{eqnarray}
We refer to this as the \textit{rotating wave solution}. An example of a rotating wave solution for $n=5$ coupled FHN oscillators in the ring configuration is presented in Figure \ref{fig:FHN_wave_ex}.

Recall equation (\ref{eq:FHN}) with $x_j=(z_j,y_j)$. If we restrict the coupled dynamics of the FHN oscillators to the rotating wave manifold, then using the periodicity of the rotating wave solution, substitution of equation (\ref{eq:rotating_wave_solution}) into the dynamics of each coupled FHN oscillator (\ref{eq:FHN}) yields $n$ identical uncoupled delay differential equations (DDEs) of the form
\begin{align*}
\dot{x}_1(t) &= f(x_1(t)) + \sigma BC\left(x_1(t-\frac{(n-1)}{n}T) - x_1(t)\right)\\
\vdots \\
\dot{x}_n(t) &= f(x_n(t)) + \sigma BC\left(x_n(t-\frac{(n-1)}{n}T) - x_n(t)\right).
\end{align*}
Thus the rotating wave solution can only exist if the auxiliary system
\begin{equation}\label{eq:aux_syst}
\begin{split}
\dot{s}(t)=&f(s(t)) -\sigma BC[s(t)-s(t-\tau^*)],\\
\tau^* :=&T-\tau=\frac{n-1}{n}T,
\end{split}
\end{equation}
has a non-constant, $T$-periodic solution:
\begin{equation}\label{eq:AuxSysEqu}
s_t=s_{t+T} \in \mathcal{C}=\mathcal{C}([0,T],\mathbb{R}^2),
\end{equation}
for which the set $\mathcal{C}$ is the set of continuous functions that map the interval $[0,T]$ into $\mathbb{R}^2$, and $s_t(\theta):=s(t+\theta)$, $\theta\in [0,T]$.

Define the errors between neighbouring systems around the rotating wave solution
\begin{equation}
\begin{array}{ll}
e_j(t) &= x_{j+1}(t+\tau) - x_j(t),\;\;j=1,2,\ldots,n-1,\\
e_n(t) &= x_1(t+\tau)-x_n(t).
\end{array}
\end{equation}
Taking the error dynamics we obtain
\begin{eqnarray}\label{eq:error_dynamics}
\left(\begin{array}{c}
\dot{e}_1(t)\\
\dot{e}_2(t)\\
\vdots\\
\dot{e}_n(t)
\end{array}\right)&=&
\left(\begin{array}{c}
f(e_1(t)+x_1(t))-f(x_1(t))\\
f(e_2(t)+x_2(t))-f(x_2(t))\\
\vdots\\
f(e_n(t)+x_n(t))-f(x_n(t))\\
\end{array}\right)-
\sigma (L\otimes BC)\left(\begin{array}{c}
e_1(t)\\
e_2(t)\\
\vdots\\
e_n(t)\end{array}\right).
\end{eqnarray}
Substitution of the rotating wave solution in terms of the auxiliary system variable $s(t)$ into equation (\ref{eq:error_dynamics}), such that
\[
s(t) = x_1(t) = x_2(t+\tau) = \cdots = x_n(t+(n-1)\tau),
\]
and linearizing around the rotating wave solution yields the linear system (\ref{eq:wide_linearized})
\begin{equation}
\left(\begin{array}{c}\label{eq:wide_linearized}
\dot{\xi}_1(t)\\
\dot{\xi}_2(t)\\
\vdots\\
\dot{\xi}_n(t)
\end{array}\right)
=
\left[\left(\begin{array}{cccc}
J(s(t)) & & & \\
 & J(s(t-\tau)) & & \\
 & & \ddots & \\
 & & & J(s(t-(n-1)\tau))\end{array}\right)
- \sigma (L\otimes BC)
\right]
\left(\begin{array}{c}
\xi_1(t)\\
\xi_2(t)\\
\vdots\\
\xi_n(t)
 \end{array}\right),
\end{equation}
where $\otimes$ is the Kronecker (tensor) product, $J(s(t))$ is defined as follows:
\[
	J(s(t)) := \begin{pmatrix}
	-\alpha \beta & \alpha \\
	-1 & 1-\gamma s_2^2(t)
	\end{pmatrix},
\]
and $s_2(t)$ denotes the second component of $s(t)$.  Note that $T$-periodicity of the system \eqref{eq:aux_syst} implies the linear error (\ref{eq:wide_linearized}) system to be $T$-periodic.

\begin{figure}
\centering
\includegraphics[width=300pt]{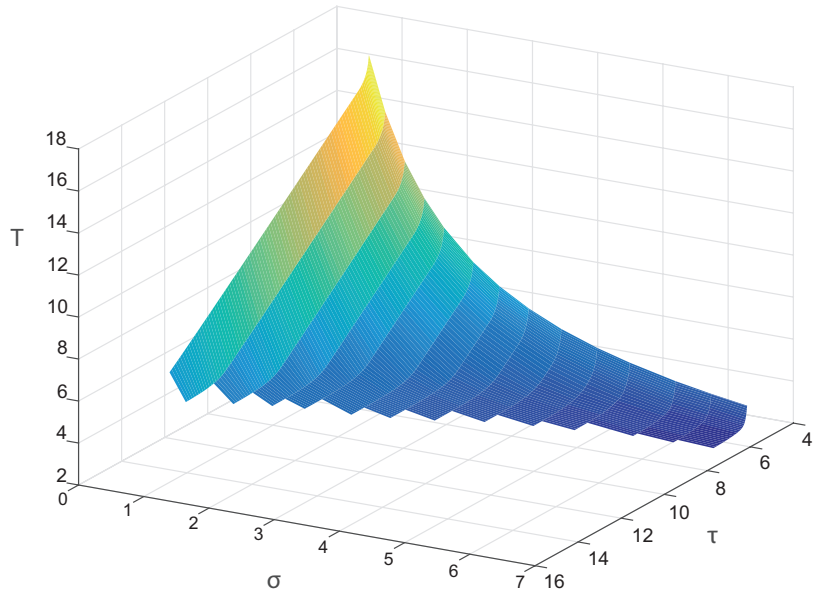}
\caption{Solutions of the auxiliary system characterized in the parameter domain of period time, delay, and coupling strength $(T,\tau,\sigma)$ presented as a surface for $T$ a function of pairs $(\tau,\sigma)$.}\label{fig:SolutionSurface}
\end{figure}

For the local stability analysis we first computed periodic solutions of the auxiliary system \eqref{eq:aux_syst}. Periodic solutions are determined using continuation methods that are available in the numerical software package \textit{DDE-Biftool} \cite{DDEBiftool}. Figure (\ref{fig:SolutionSurface}) characterizes solutions of the auxiliary system in the parameter domain  $(T,\tau,\sigma)$.  For the auxiliary system in which parameter $T$ and $\tau$ are allowed to vary continuously, a solution that describes the dynamics of a rotating wave solution satisfies the relation $\frac{T}{\tau}(n-1) = n$. However, for the solutions we obtained, parameters $T$ and $\tau$ have not been varied continuously. Therefore, we choose the solution that satisfies the following the inequality
\begin{equation}\label{eq:Ttau_to_n}
\left| \frac{T}{\tau}(n-1) - n \right| < \epsilon.
\end{equation}
To maintain good accuracy of approximation of the auxiliary system to $n$ coupled FHN oscillators, the error $\epsilon$ must be small. For our stability analysis we took $\epsilon = 0.01$.

\begin{figure}[!h]
\subfigure[]{\label{fig:sigma095}
\includegraphics[width=0.45\textwidth]{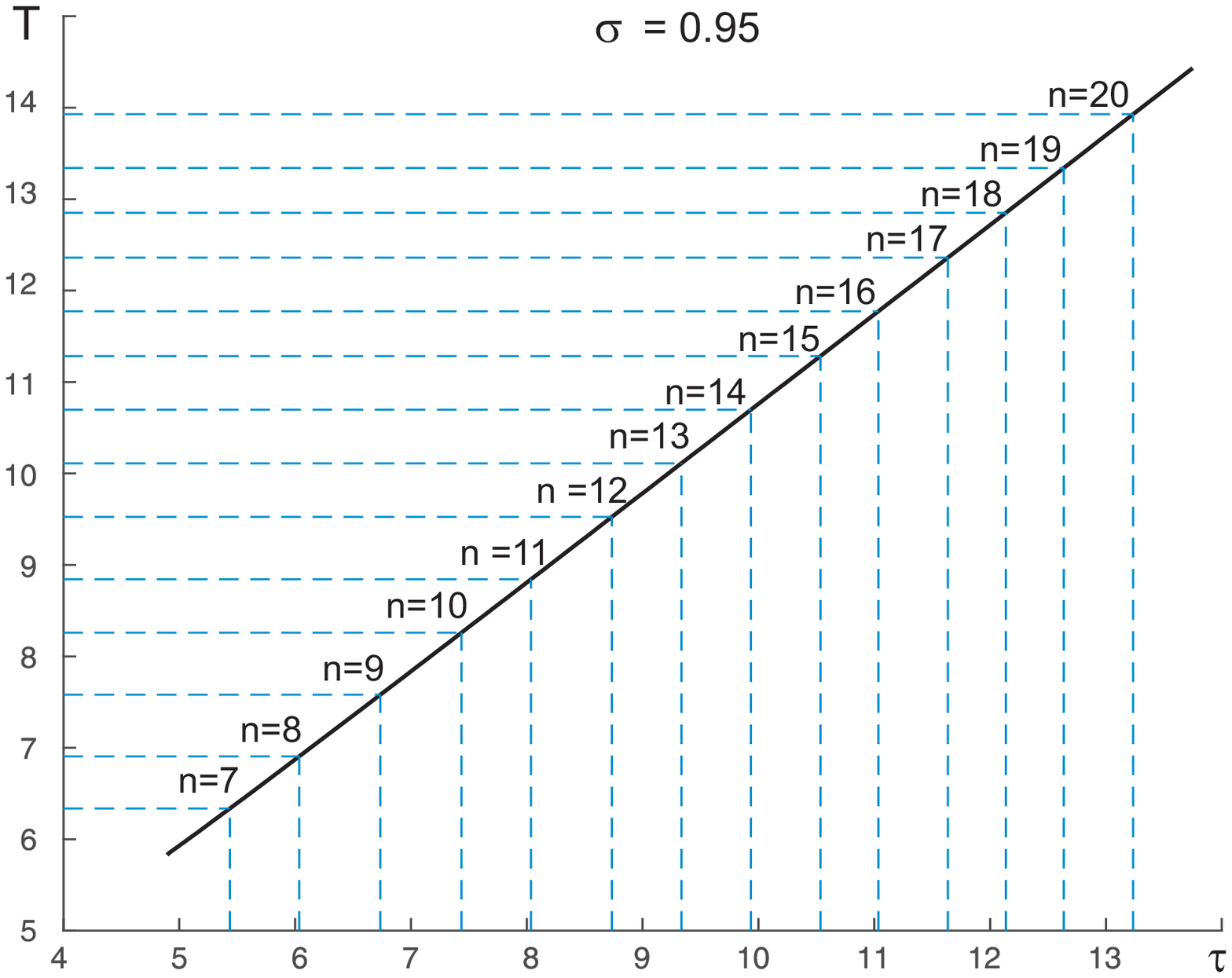}
}\hfill
\subfigure[]{\label{fig:sigma675}
\includegraphics[width=0.45\textwidth]{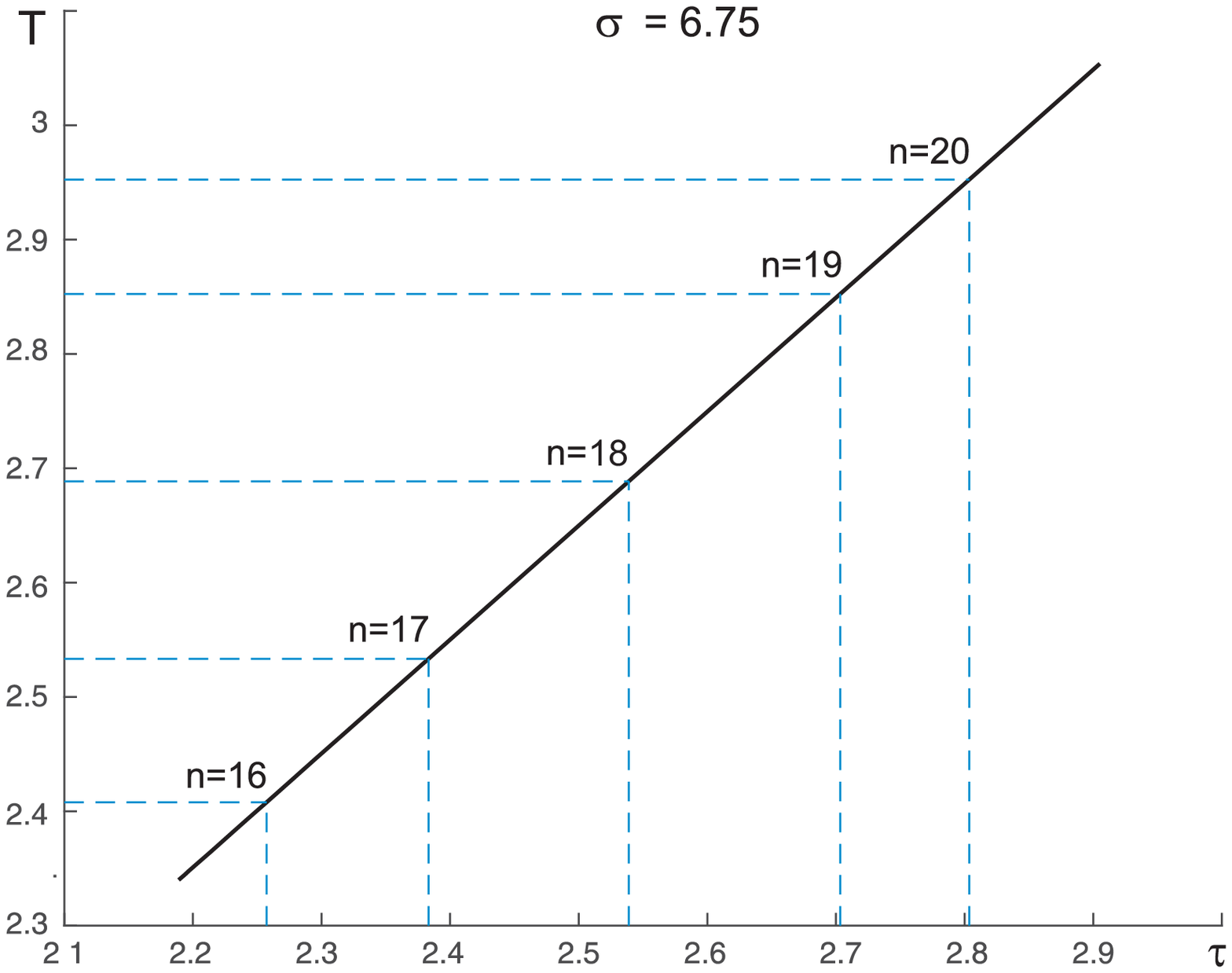}
}
\caption{Solutions of the auxiliary system in the parameter domain of period time and delay $(T,\tau)$ for given coupling strengths: (a) $\sigma=0.95$. (b) $\sigma=6.75$. Dashed lines indicate solutions that satisfy relation (\ref{eq:Ttau_to_n}).}\label{fig:T_over_tau}
\end{figure}

Figures \ref{fig:sigma095} and \ref{fig:sigma675} show two cross sections of the surface in Figure \ref{fig:SolutionSurface} for coupling strengths $\sigma=0.95$ and $\sigma=6.75$, respectively. Dashed lines identify solutions that satisfy relation (\ref{eq:Ttau_to_n}) and hence map solutions of the auxiliary system to an integer number $n$ of coupled FHN oscillators on the rotating wave manifold.

We assessed the stability of the rotating wave solution for pairs $(n,\sigma)$ by computing the Floquet multipliers of the periodic linearized error system (\ref{eq:wide_linearized}) (again with \textit{DDE-Biftool}). Solutions are obtained by substituting in the solution of the auxiliary system corresponding to the pair $(n,\sigma)$ found by numerical continuation. Recall that if all Floquet multipliers except one (at $1$) have modulus strictly smaller than $1$, then the zero solution of the linearized error system is asymptotically stable, which implies the rotating wave solution to be locally orbitally stable. The red line in Figure \ref{Fig:FHN_network} (and Figure \ref{Fig:FHN_network_wave_basin_type1}) is defined by the crossing of (at least) one multiplier with the boundary of the unit disc in $\mathbb{C}$.

\setcounter{equation}{0}

\section{Discussion}\label{sec:discussion}

\subsection{Kinetic interpretation of (\ref{eq:sysP})}

Equation \eqref{eq:sysP} in Section \ref{sec:eigenvalues} describes the temporal evolution of the first order kinetics. This equation is known as the {\em Master Equation}.
The master equation obeys the principle of detailed balance if there exists a positive
equilibrium $P^*$ ($p^*_i>0$) such that for each pair $i,j$ ($i\neq j$)
\begin{equation}\label{detbal}
q_{ij} p^*_j=q_{ji} p^*_i .
\end{equation}
After Onsager \cite{Ons}, it is well known that for systems with detailed balance the
eigenvalues of $K$ are real because under conditions (\ref{detbal}) $K$ is a
self-adjoined matrix with respect to the entropic inner product
$$\langle x,y\rangle=\sum_i \frac{x_i y_i}{p^*_i}$$
(see, for example, \cite{VanKampen1973,YBGE1991}).

Detailed balance is a well known consequence of microreversibility. This principle was introduced in 1872
by Boltzmann for collisions \cite{Boltzmann1964}. In 1901 Wegscheider proposed
it for chemical kinetics \cite{Wegscheider1901}. Einstein had used it as a principle
for the quantum theory of light emission and absorption (1916, 1917). The
backgrounds of detailed balance had been analyzed by Tolman \cite{Tolman1938}. The
principle was studied further and generalized by several authors
\cite{Gorban2014,YangHlavacek2006,GorbYabCES2012}.

Systems without detailed balance appear in applications rather often. Usually, they
represent a subsystem of a larger system, where concentrations of some of the components
are considered as constant. For example, the simple cycle
\begin{equation}\label{simplecycle}
A_1 \to A_2 \to \ldots \to A_n \to A_1
\end{equation}
is a typical subsystem of a catalytic reaction (a catalytic cycle). The complete reaction
may have the form
\begin{equation}\label{Catcycle}
S+A_1 \to A_2 \to \ldots \to A_n \to A_1+P,
\end{equation}
where $S$ is a substrate and $P$ is a product of reaction.

The irreversible cycle (\ref{simplecycle}) cannot appear as a limit of systems with
detailed balance when some of the constants tend to zero, whereas the whole catalytic
reaction (\ref{Catcycle}) can \cite{GorbYabCES2012}. The simple cycle (\ref{simplecycle})
can be produced from the whole reaction (\ref{Catcycle}) if we assume that concentrations
of $S$ and $P$ are constant. This is possible in an open system, where we continually add
the substrate and remove the product. Another situation when such an approximation makes
sense is a significant excess of substrate in the system, $[S]\gg [A_i]$ (here we use the
square brackets for the amount of the component in the system). Such excess implies
separation of time and the system of intermediates $\{A_i\}$ relaxes much faster than the
concentration of substrate changes.

In systems without detailed balance, damped oscillations are possible. The example in Section \ref{sec:eigenvalues}, which describes the case of all the reaction rate constants in the simple cycle being equal, $q_{j+1\, j}=q_{1n}=q>0$, shows that these oscillations are even weakly damped. The effect becomes acutely noticeable for $n$ large enough.


The simple cycle with equal rate constants yields the slowest decay of oscillations or, in some sense, the slowest relaxation
among all first order kinetic systems with the same number of components. The extremal properties of the simple cycle with equal constants
were noticed in numerical experiments $25$ years ago \cite{BochByk1987}. V.I. Bykov formulated the hypothesis that this
system has extremal spectral properties. This paper provides the answer: yes, it has.

%

\subsection{Two coupled cycles}

Given the size of the region where multiple solutions co-exist, and the resilience to a coherent state; does the extremal property of the simple cycle give rise to further, more complex phenomena when two simple cycles are diffusively coupled via an undirected link between an oscillator in each cycle?

For a total of $2k$ coupled systems, two cycles are constructed with systems $1,\ldots,k$ in the first simple cycle and systems $k+1,\ldots,2k$ in the second, and coupled via systems $x_1$ and $x_{k+1}$. Clearly the synchronization manifold exists, as does the rotating wave solution in the form of two synchronized rotating waves,
\begin{multline*}
x_1(t) = x_{k+1}(t) = x_2(t+\tau) = x_{k+2}(t+\tau) = \ldots \\
\ldots = x_{k-1}(t+(k-2)\tau) = x_{2k-1}(t+(k-2)\tau) = x_k(t+(k-1)\tau) = x_{2k}(t+(k-1)\tau).
\end{multline*}
A full description of the phenomena of two coupled cycles is beyond the scope of this work; however, as a motivation for further study, we present a brief example.

\begin{figure}
\centering{
\includegraphics[width=250pt]{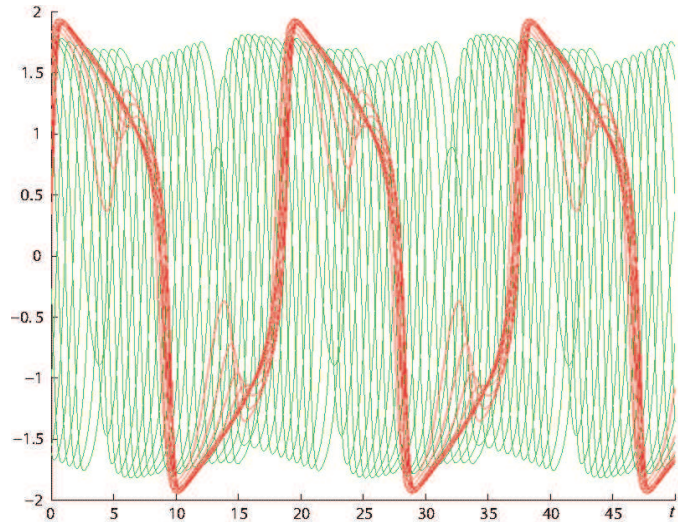}
\caption{\label{Fig:sync_and_wave} Two coupled cycles and their \textit{y}-dynamics; in red the \textit{y}-dynamics of the first cycle, and in green the \textit{y}-dynamics of the second cycle.}
}
\end{figure}

We take $(n,\sigma)=(10,0.75)$, which, for a simple cycle lies in the  region of co-existence of synchronization and rotating wave solutions. We observe in Figure (\ref{Fig:sync_and_wave}) a stable state in which the trajectories of all systems in the first cycle (in red) are attracted to the synchronization manifold, whilst all trajectories of systems in the second cycle (in green) are attracted to the rotating wave solution. There is a clear competition of each cycle to attract the other to its own dynamical regime.
 The two diffusively coupled oscillators from each cycle periodically perturb each other, which prevents asymptotic convergence of systems to either the synchronization manifold or the rotating wave solution. Clearly, the extremal properties of the simple cycle can give rise to multiple regimes of complex patterns of dynamics when embedded into larger network structures.

\section{Conclusion}\label{sec:conclusion}

We considered the problem of how ``closing'' a chain of interconnected systems with directed coupling by adding a directed feedback from the last element in the chain to the first may affect the dynamics of the system. This problem is closely related to the fundamental question of how network topology influences the dynamics of collective behavior in the system. Two general settings have been investigated. In the first one we analyzed the behavior of a simple linear system. We showed that the simple cycle with equal interaction weights has the slowest decay of the oscillations among all linear systems with the same number of states. In the second setting we considered directed rings and chains of identical nonlinear oscillators. For  directed rings, a lower bound $\sigma_c$ for the connection strengths that guarantee asymptotic synchronization in the network is found to follow a pattern similar to that of a simple cycle. Furthermore, numerical analysis revealed that, depending on the network size $n$, multiple dynamic regimes co-exist in the system's state space.

In addition to the fully synchronous state, for sufficiently large networks an asymptotically stable rotating wave solution emerges. The emergence of the rotating wave is a phenomenon that persists over a broad range of coupling strengths and network sizes, and can be viewed as a form of extreme sensitivity of the network dynamics to the removal or addition of a single connection. The result confirms the significance of shortcuts in  networks with large numbers of nodes. Emergence of asymptotically stable rotative wave solutions has been analyzed numerically for a specific class of systems in which the dynamics of each node was identical and satisfied Fitzhugh-Nagumo equations \cite{FitzHugh}. Extending the analysis to systems with heterogeneous nodes as well as considering nodes with Hindmarsh-Rose and Hodgkin-Huxley dynamics \cite{Izhikevich}, known to be capable of bursting and chaotic behavior, will be the topic of our future studies.

Coming back to the question if leaders should look back. To stay in synchrony we advise a leader either not to look back at all or to look back just a few links; looking back too far induces oscillations that destroy the coherent state.

\begin{acknowledgement}
The authors are thankful to anonymous Referees for their encouraging and helpful suggestions and comments. Ivan Tyukin is also thankful to the Russian Foundation for Basic Research (research project No. 15-38-20178) for partial support. Cees van Leeuwen was supported by an Odysseus Grant from the Belgion Foundation for Science, F.W.O.
\end{acknowledgement}

\end{document}